\documentclass[12pt]{article}
\usepackage{amsmath,amsthm,amssymb,latexsym, amsfonts}
\usepackage[dvipsnames]{xcolor}
\usepackage[dvips]{graphicx}

\theoremstyle{plain}
\newtheorem{teor}{Theorem}[section]
\newtheorem{defin}[teor]{Definition}
\newtheorem{propo}[teor]{Proposition} 
\newtheorem{obs}[teor]{Remark} 
\newtheorem{theo}[teor]{Theorem}
 
\newtheorem{example}[teor]{Example} 
\newtheorem{coro}[teor]{Corollary} 

\def\proof{{\noindent \bf Proof:} \hspace{0.1 cm}}
\newcommand{\cqd}{\hfill \rule{2mm}{2mm}\vspace{.3cm}}

\date{}
\title{About quadratic residues in a class of rings}

\begin{document}

\author{ Fernanda D. de Melo Hern\'andez \thanks{Corresponding author}
\thanks{F. D. de Melo Hernandez (fdmelo@uem.br), Departamento de Matem\'atica, Universidade Estadual de Maring\'a, Av. Colombo 5790, 87020-900, Maring\'a, PR, Brazil.} \and
C\'esar A. Hern\'andez Melo \thanks{C\'esar A. Hern\'andez Melo (cahmelo@uem.br), Departamento de Matem\'atica, Universidade Estadual de Maring\'a, Av. Colombo 5790, 87020-900, Maring\'a, PR, Brazil.} 
\and Horacio Tapia-Recillas \thanks{Horacio Tapia-Recillas (htr@xanum.uam.mx), Departamento de Matem\' aticas, Universidad Aut\' onoma Metropolitana-Iztapalapa, Av. San Rafael Atlixco 186, 09340, CDMX, M\' exico.}}
\maketitle


\begin{abstract}

Let $R$ be a commutative ring with a collection of ideals $\{ N_1, N_2, \dots, N_{k-1}\}$ satisfying certain conditions, properties of the set of invertible quadratic residues of the ring $R$ are described in terms of properties of the set of invertible quadratic residues of the quotient ring $R/N_1$.


\end{abstract}



\section{Introduction}

Quadratic residues (\cite{i-r}) have been studied since the 17th and 18th century by P. de Fermat, L. Euler, J.L. Lagrange A.M. Legendre, among other mathematicians. Nowadays quadratic residues are still a topic of study (\cite{c-f}),\cite{l-s},\cite{b1},\cite{p-s},\cite{s}) and they have applications in areas which include acoustical engineering (\cite{sch}), cryptography (\cite{j},\cite{h-k}), in the study of Paley (conference) graphs, primality testing (Solovay-Strassen, Miller-Rabin), and integer factorization (quadratic sieve, number field sieve). \\

In the present note properties of invertible quadratic residues over a class of rings are given 
 which are determined by properties of quadratic residues over a quotient ring, i.e., the properties are ``lifted" from those of a quotient ring. An example of this class of rings includes the integers modulo $p^k$, $\mathbb Z_{p^k}$, where $p$ is a prime and $k$ a positive integer.\\
 
 The manuscript is divided in four sections. In Section 2 notation and facts needed in the rest of the manuscript are presented. In section 3 the main results are given and in Section 4 applications of the results previously discussed are considered. Examples are given to illustrate the main results.



\section{Notations and basic facts}\label{BF}

Given $R$ an associative ring with identity and $N$ a nil ideal of $R$, we begin our discussion  recalling that units of the quotient ring $R/N$ can be lifted to the ring $R$. More precisely, if $R^*$ and $(R/N)^*$ denote the group of units of the ring $R$ and $R/N$ respectively, the following result holds. 
\begin{propo}\label{jacobson}
Let $R$ be an associative ring with identity, $N$ a nil ideal of $R$ and  \hspace{2mm}  $\bar{} : R \longrightarrow R/N$ the canonical homomorphism from $R$ to the quotient ring $R/N$. Then, 
\begin{enumerate}
\item $\bar f=f+N\in (R/N)^{\ast}$ if and only if $f+N\subset R^{\ast}$.  
\item If $R$ is finite the cardinality of $R^*$ and the cardinality of $(R/N)^*$ are related by the relation
\begin{equation}\label{fcunits}
|R^*|=|N||(R/N)^*|.
\end{equation}
\end{enumerate}
\end{propo}
\proof The proof of this proposition can be found in \cite{liftunits}, Proposition 2.1 and Remark 2.2.
\cqd

Recall that an element $a$ of a ring  $R$ is a \textit{quadratic residue}, if there exists $x\in R$ such that $x^2=a$ (\cite{i-r},\cite {n-z},\cite{w1}). Given $N$ an ideal of $R$, it is clear that if $a$ is a quadratic residue in the ring $R$, then $a+N$ is a quadratic residue in the quotient ring $R/N$. The following proposition provides sufficient conditions to prove the converse of this statement, as it will be seen, proposition \ref{jacobson} will be essential in the proof that will be presented.
\begin{propo}\label{quadN2}
Let $R$ be a commutative ring with identity and $N$ a nil ideal of $R$. If $(g+N)^2=a+N$ and $2g+N$ is an invertible element in $R/N$, then the function 
$$
\eta:g+N\rightarrow a+N\hspace{0.3cm}\text{given by}\hspace{0.3cm}\eta(g+m)=(g+m)^2,
$$ is bijective. In other words, if $a+N$ is a quadratic residue in $R/N$, then every element $b\in a+N$ is a quadratic residue in the ring $R$ and, and for all $b \in a+N$ the quadratic equation
$$
y^2=b
$$
has a unique solution in the set $g+N\subset R$.
 
\end{propo}
\proof 
Since $(g+N)^2=a+N$, it is clear that the function $\eta$ is well defined. Since $N$ is a nil ideal of $R$ and $2g+N$ is an invertible element in $R/N$, from proposition \ref{jacobson}, it follows that for all $p\in N$, the element $2g+p$ is an invertible element in $R$. Thus, if $\eta(g+m_1)=\eta(g+m_2)$ then, 
$$
(2g+m_1+m_2)(m_1-m_2)=0.
$$ 
Since $2g+m_1+m_2$ is an invertible element in the ring $R$, it is concluded that $m_1=m_2$. Therefore $\eta$ is an injective function. Now, it is proved that $\eta$ is surjective. First of all, note that since $(g+N)^2=a+N$, there exists $n_0\in N$, such that $g^2=a+n_0$. Now, given $n\in N$, it is easy to see that    
$$
\eta(g+(2g)^{-1}(n-n_0))=a+n,
$$
which proves the claim, finishing the proof of the proposition.
\cqd
\\
Now, definitions and notation that will be useful in the rest of the manuscript are introduced. The set $q(R^*)$ will denote the quadratic residues in the ring $R$ that are also units in $R$, that is 
$$
q(R^*)=\{a\in R| a\hspace{0.2cm}\text{is a quadratic residue in}\hspace{0.2cm}R\hspace{0.2cm}\text{and}\hspace{0.2cm} a\in R^* \}.
$$
For $a$ a quadratic residue in the ring $R$, $s(a)$ will denote the set of solutions of the equation $x^2=a$ in the ring $R$. In other words,
$$
s(a)=\{x\in R| x^2=a\}.  
$$ 
Finally, if $N$ is an ideal of the ring $R$ and $a\in R$, $T(a+N)$ will be denote the set of solutions of the equation $x^2=b$, when $b$ varies in the equivalence class $a+N\in R/N$, in other words
$$
T(a+N)=\{y\in R| y^2\in a+N\}.
$$
Based on propositions \ref{jacobson} and \ref{quadN2}, sufficient conditions to lift quadratic residues from the quotient ring $R/N$ to the ring $R$, where $N$ is a nil ideal of the ring $R$ are established. In addition if $R$ is finite, formulas relating the cardinality of the sets $N, s(b), s(b+N), R^*, (R/N)^*, q(R^*)$ and $q((R/N)^*)$ are given. 

\begin{propo}\label{coroe}
Let $R$ be a commutative ring with identity and $N$ a nil ideal of $R$ such that $2+N$ is an invertible element in $R/N$. The following statements hold,
\begin{enumerate}
\item\label{ide1} $a+N\in q((R/N)^*)$ if and only if $a+N\subset q(R^*)$.
\item\label{ide3} The cardinality of the set $q(R^*)$ satisfies
\begin{equation}\label{ces2}
|q(R^*)|=|N||q((R/N)^*)|.
\end{equation}
\item\label{ide2} If $a+N\in q((R/N)^*)$, then for all $b\in a+N$ the number of solutions of the quadratic equation 
$$
y^2=b
$$ in the ring $R$ is equal to the number of solutions of the quadratic equation 
$$
y^2=b+N
$$ in the ring $R/N$. In other words 
$$
|s(b)|=|s(a+N)|
$$
for all $b\in a+N$.
\item\label{ide3.1} The cardinality of the set $R^*$ satisfies the following relation
\begin{equation}\label{ces1}
|R^*|=|N|\sum_{a+N\in q((R/N)^*)} |s(a+N)|.
\end{equation}
\item\label{ide5} If in addition, there exists $\alpha$ such that $|s(a+N)|=\alpha$ for all $a+N\in q((R/N)^*)$,
\begin{equation}\label{ces3}
\it{a})\hspace{0.1cm}|q((R/N)^*)|=\frac{|(R/N)^*|}{\alpha}, \hspace{0.6cm}\it{b})\hspace{0.1cm}|q(R^*)|=\frac{|N||(R/N)^*|}{\alpha},\hspace{0.6cm}\it{c})\hspace{0.1cm}|q(R^*)|=\frac{|R^*|}{\alpha}.
\end{equation}

\end{enumerate}   
\end{propo} 
\proof \textit{\ref{ide1}}. It is easy to see that if $a+N\subset q(R^*)$, then $a+N\in q((R/N)^*)$. Now, it is proved the other implication. Assuming $a+N\in q((R/N)^*)$, there exists $g+N\in R/N$ such that $(g+N)^2=a+N$. Since $a+N$ and $2+N$ are invertible elements in the ring $R/N$, it follows that
$$
(2+N)(g+N)=2g+N
$$ is an invertible element in $R/N$, thus proposition \ref{quadN2} implies that every element $b\in a+N$ is a quadratic residue in the ring $R$. In addition, since $a+N\in (R/N)^*$ and $N$ is a nil ideal of the ring $R$, proposition \ref{jacobson} implies that for all $b\in a+N$, $b\in R^*$. This proofs that $a+N\subset q(R^*)$.\\
\textit{\ref{ide3}}. From claim \textit{\ref{ide1}}, it follows that
\begin{equation}\label{cardrq}
q(R^*)=\bigcup_{a+N\in q((R/N)^*)} (a+N).
\end{equation}
Then,
$$
|q(R^*)|=\sum_{a+N\in q((R/N)^*)} |a+N|=|N|\sum_{a+N\in q((R/N)^*)} 1=|N||q((R/N)^*)|,
$$
which proves relation (\ref{ces2}).
\\
\textit{\ref{ide2}}. Note that since $a+N\in q((R/N)^*)$, for all $b\in a+N$, $s(b)\neq\emptyset$. In order to prove that $|s(b)|=|s(a+N)|$, it will be shown that the canonical homomorphism $\phi: R\rightarrow R/N$ restricted to the set $s(b)$, namely
$$
x\in s(b)\rightarrow \phi(x)=x+N
$$
determines a bijection between $s(b)$ and $s(b+N)$. In fact, if $x\in s(b)$ then $x+N\in s(b+N)$, thus the function $\phi$ is well defined. Now, if $\phi(x)=\phi(y)=x+N$ with $x,y\in s(b)$, then $(x+N)^2=b+N=a+N$  with $2x+N$ an invertible element in the ring $R/N$. So, proposition \ref{quadN2} implies that the function 
$$
z\in x+N\rightarrow \eta_1(z)=z^2\in b+N
$$
is a bijective function, hence, since $x,y\in x+N$ and $x^2=y^2=b$, the injectivity of the function $\eta_1$ implies that $x=y$. Now, if $t+N\in s(b+N)$, then $(t+N)^2=b+N=a+N$. Since $2t+N$ is an invertible element in the ring $R/N$, proposition \ref{quadN2} implies that the function 
$$
z\in t+N\rightarrow \eta_2(z)=z^2\in b+N
$$
is bijective. Thus there exists $n_0\in N$, such that $\eta_2(t+n_0)=(t+n_0)^2=b$. Hence, $t+n_0\in s(b)$ and $\phi(t+n_0)=t+N$, hence $\phi$ is a surjective function.  
 \\
\textit{\ref{ide3.1}}. Note that the set $R^*$ is a disjoint union of the sets $T(a+N)$ with $a+N\in q((R/N)^*)$, that is
$$
R^*=\bigcup_{a+N\in q((R/N)^*)} T(a+N).
$$
From this fact, it follows that
\begin{equation}\label{gew}
|R^*|=\sum_{a+N\in q((R/N)^*)} |T(a+N)|.
\end{equation}
In order to compute $|T(a+N)|$, observe that $T(a+N)$ can be written as a disjoint union of the sets $s(b)$ with $b\in a+N$, 
\begin{equation}\label{cardrqsol}
T(a+N)=\bigcup_{b\in a+N} s(b).
\end{equation}
Thus, since for all $b\in a+N$, $|s(b)|=|s(a+N)|$, 
\begin{equation}\label{gew1}
|T(a+N)|=\sum_{b\in a+N} |s(b)|=|s(a+N)|\sum_{b\in a+N} 1=|s(a+N)||N|.
\end{equation}
Finally, combining (\ref{gew}) and (\ref{gew1}), relation in (\ref{ces1}) follows easily.
 \\
\textit{\ref{ide5}}. Since $|s(a+N)|=\alpha$ for all $a+N\in q((R/N)^*)$, from relation (\ref{ces1}), it follows that
\begin{equation}\label{pri1}
|R^*|=\alpha|N||q((R/N)^*)|,
\end{equation}
and proposition (\ref{jacobson}) implies that, 
\begin{equation}\label{pri2}
|R^*|=|N||(R/N)^*|.
\end{equation}
Thus, by combining (\ref{pri1}) and (\ref{pri2}), relation in (\ref{ide5})-\textit{a}) is obtained. The relation in (\ref{ide5})-\textit{b}) is obtained from (\ref{ces2}) and (\ref{ide5})-\textit{a}). Finally, relation in (\ref{ide5})-\textit{c}) is obtained from (\ref{pri2}) and (\ref{ide5})-\textit{b}).
\cqd
\\
In the next lines the results in proposition \ref{coroe} are illustrated with an example. Let $p$ be a prime number different from $2$, $k$ a natural number and let $R=\mathbb{Z}_{p^k}$ be the ring of integers modulo $p^k$. It is clear that the ideal generated by $p$ in $R$, namely $N=\langle p \rangle$, is a nilpotent ideal of index $k$ of the ring $R$. In addition, 
$$
\frac{R}{N}=\frac{\mathbb{Z}_{p^k}}{\langle p \rangle}\cong \mathbb{Z}_p,
$$
thus $|(R/N)^*|=p-1$ and from Lagrange's theorem it follows that $|N|=p^{k-1}$. From the identity in (\ref{fcunits}), we have that $|R^*|=p^{k-1}(p-1)$. In addition, since, $R/N\cong \mathbb{Z}_p$ is a field of characteristic different from 2, it follows that the number $\alpha$ appearing in claim \textit{\ref{ide5}} of proposition \ref{coroe} is $\alpha=2$. Thus, from proposition \ref{productring}, it is concluded that:
\begin{itemize}
\item If $a\in\mathbb{Z}_{p^k}^*$, and $a\equiv b\mod(p)$, then $a$ is a quadratic residue in $\mathbb{Z}_{p^k}$ if and only if $b$ is a quadratic residue in the ring $\mathbb{Z}_{p}$.
\item If $a\in q(\mathbb{Z}_{p^k}^*)$ and $a\equiv b\mod(p)$, then the number of solutions of the equation $x^2=a$ in the ring $\mathbb{Z}_{p^k}$ is equal to the number of solutions of the equation $x^2=b$ in the field $\mathbb{Z}_p$, which is equal to $2$, in other words
$$
s(a)=s(b)=2.$$
\item The cardinality of the sets $q(\mathbb{Z}_p^*)$, $q(\mathbb{Z}_{p^k}^*)$ are given by
$$
|q(\mathbb{Z}_p^*)|=\frac{p-1}{2}\hspace{0.6cm}\text{and}\hspace{0.6cm}|q(\mathbb{Z}_{p^{k}}^*)|=\frac{p^{k-1}(p-1)}{2},
$$
 respectively.
\end{itemize}

Next, the previous proposition is extended to a direct product of a finite collection of rings.
\begin{propo}\label{productring} 
Let $R_1, R_2,\dots,R_m,$ be commutative rings with identity and let $N_i$ be a nil ideal of the ring $R_i$, such that $2+N_i\in (R_i/N_i)^*$ for each $i=1,2,\dots, m$. The following statements hold:
\begin{enumerate}
\item\label{idea1} $(a_1,\dots,a_m)\in q((R_1\times\cdots\times R_m)^*)$ if and only if $a_i+N_i\in q((R/N_i)^*)$ for every i=1,2,\dots,m.
\item\label{idea2} If $(a_1,\dots,a_m)\in q((R_1\times\cdots\times R_m)^*)$ then 
\begin{equation}\label{mist}
|s((a_1,\dots,a_m))|=|s(a_1+N_1)|\cdots |s(a_m+N_m)|.
\end{equation}
\item\label{idea3} The cardinality of $q((R_1\times\cdots\times R_m)^*)$ satisfies the following relation
\begin{equation}\label{tela}
|q((R_1\times\cdots\times R_m)^*)|=|N_1||q((R/N_1)^*)|\cdots|N_m||q((R/N_m)^*)|
\end{equation}
\item\label{idea4} If $|s(a+N_i)|=\alpha_i$ for all $a+N_i\in q((R/N_i)^*)$, then
\begin{equation}\label{ces31}
|q((R_1\times\cdots\times R_m)^*)|=\frac{|N_1||(R/N_1)^*|\cdots |N_m||(R/N_m)^*|}{\alpha_1\alpha_2\cdots\alpha_m},
\end{equation}
and
\begin{equation}\label{ces41}
|q((R_1\times\cdots\times R_m)^*)|=\frac{|(R_1\times\cdots\times R_m)^*|}{\alpha_1\alpha_2\cdots\alpha_m}.
\end{equation}
\end{enumerate}  
\end{propo}
\proof \textit{\ref{idea1}}. The proof of this claim follows from claim \textit{\ref{ide1}}) of proposition \ref{coroe} and the fact that
\begin{equation}\label{qprod}
q((R_1\times\cdots\times R_m)^*)=q(R_1^*)\times\cdots\times q(R_m^*).
\end{equation}
\textit{\ref{idea2}}. Since $s((a_1,\dots,a_m))=s(a_1)\times\cdots\times s(a_m)$,
$$
|s((a_1,\dots,a_m))|=|s(a_1)|\cdots|s(a_m)|.
$$
Thus, relation (\ref{mist}) follows from \textit{\ref{idea1}}) of proposition \ref{productring} and claim \textit{\ref{ide2}}) of proposition \ref{coroe}.
\\
\textit{\ref{idea3}}. The proof of relation in (\ref{tela}) is a consequence of the equality given in (\ref{qprod}) and the identity in (\ref{ces2}).\\
\textit{\ref{idea4}}. Relation (\ref{ces31}) follows from (\ref{tela}) and claim (\ref{ide5})-\textit{a}) of proposition \ref{coroe}. Finally, since $(R_1\times\cdots\times R_m)^*=R_1^*\times\cdots\times R_m^*$, 
relation (\ref{ces41}) follows from the fact that $|R_i^*|=|N_i||(R/N_i)^*|$ for every $i=1,2,\dots,m$.
\cqd 

In the following lines the results of the previous proposition are illustrated with an example. Given $n$ an odd natural number,  if $n=p_1^{k_1}p_2^{k_2}\cdots p_m^{k_m}$ denotes the prime factorization of $n$. The Chinese Remainder Theorem implies that
$$
\mathbb{Z}_n\cong \mathbb{Z}_{p_1^{k_1}}\times\cdots\times\mathbb{Z}_{p_m^{k_m}}.
$$
By setting $R_i=\mathbb{Z}_{p_i^{k_i}}$ and $N_i=\langle p_i\rangle$ for $i=1,2,\dots,m$, it is clear that $R_i/N_i\cong \mathbb{Z}_{p_i}$ and that $2+N_i\in (R_i/N_i)^*$. Thus, from proposition \ref{productring}, it is concluded that:
\begin{itemize}
\item If $a\in\mathbb{Z}_n^*$, and $a\equiv a_i\mod(p_i)$ for $i=1,2,\dots,m$, $a$ is a quadratic residue in $\mathbb{Z}_n$ if and only if $a_i$ is a quadratic residue in the ring $\mathbb{Z}_{p_i}$ for all $i=1,2,\dots,m.$
\item If $a\in q(\mathbb{Z}_{n}^*)$ and $a\equiv b_i\mod(p_i^{k_i})$ for $i=1,2,\dots,m$, since $s(b_i+N_i)=2$ for all $i=1,2,\dots,m$, the number of solutions of the equation $x^2=a$ in the ring $\mathbb{Z}_{n}$ is equal to $2^m$, in other words
$$
s(a)=2^m.$$
\item The cardinality of the set $q(\mathbb{Z}_n^*)$ is given by
$$
|q(\mathbb{Z}_n^*)|=\frac{p_1^{k_1-1}(p_1-1)\cdots p_m^{k_m-1}(p_m-1)}{2^m}.
$$
\end{itemize}


\section{Main results}

In this section the main results of the manuscript are presented. For $R$ a commutative ring containing a collection of ideals $\{ N_1, N_2, \dots, N_{k-1}\}$ satisfying a certain condition (the CNC condition, Definition \ref{PosLiftIde}), properties of the set of invertible quadratic residues of the ring $R$ are described in terms of properties of the set of invertible quadratic residues of the quotient ring $R/N_1$.\\

\begin{propo} \label{potencia}
Let $R$ be a commutative ring and $N$ a nilpotent ideal of index $t \geq 2$ in $R$. Then the following statements hold:
\begin{enumerate}
	\item\label{pot} For any prime number $p$ such that $p \geq t$, for all $n\in N$ and $a\in R$,  
	$$(a+n)^p = a^p + pnr,$$  
	for some $r\in R$.
	\item\label{port} In addition, assuming there exists a natural number $s>1$ such that $sN=\{0\}$ and such that all the prime factors of $s$ are greater than or equal to the nilpotency index $t$ of the ideal $N$. The function $H:R/N\rightarrow R$ given by 
	$$
	H(x+N)=x^s
	$$ is well defined  and it is multiplicative, i.e., it satisfies $H((x+N)(y+N))=H(x+N)H(y+N),$ for all $x,y\in R$.  
		\item\label{porta} Under the assumptions of claim \ref{port}, if $a+N$ is a quadratic residue in the quotient   
	ring $R/N$, then $H(a+N)=a^s$ is a quadratic residue in $R$. More precisely, if $g\in R$ is such that $(g+N)^2=a+N$, then
	\begin{equation*}\label{quadraticlema}
	\left(g^s\right)^2=a^s.
	\end{equation*}
\end{enumerate}
\end{propo}

\proof 
\begin{enumerate}
\item  Since $n^t = 0$,
$$(a + n)^p = \sum_{j=0}^{p}\binom{p}{j}a^{p-j}n^j = a^p + \sum_{j=1}^{t-1}\binom{p}{j}a^{p-j}n^j.$$
Since $p$ is a prime number, $p$ divides $\binom{p}{j}$ for all $1\leq j\leq p-1$.  Also, since $t\leq p$,
$$(a+n)^p=a^p+pn\left(k_1a^{p-1}+{k_2}a^{p-2}n+\cdots + k_{t-1}a^{p-t+1}n^{t-2}\right)$$
where $k_i=\tbinom{p}{i}/p$. Therefore,
$$(a+n)^{p}= a^p+pnr,$$
with $r=k_1a^{p-1}+k_2a^{p-2}n+\cdots + k_{t-1}a^{p-t+1}n^{t-2}\in R$.

\item Let $p_1,p_2,p_3,\ldots ,p_m$ be the prime numbers, not necessarily different, appearing in the prime factor decomposition of the integer $s$, with $p_i \geq t$, for $i=1,2,3,\cdots,m.$ Since $y+N=x+N$, there exists $n \in N$ such that $y = x + n$. Since $p_1\geq t$, from claim \ref{pot}, 
$$y^{p_1}=(x + n)^{p_1}= x^{p_1} + p_1nr_1,$$
for some $r_1\in R$. Similarly, since $p_2\geq t$ and $p_1nr_1\in N$, it follows from claim \ref{pot} and the previous relation that  
$$y^{p_1p_2}=(x^{p_1}+p_1nr_1)^{p_2}=x^{p_1p_2}+p_2(p_1nr_1)r_2,$$
for some $r_2\in R$. In the same way, it is possible to verify that  
$$y^{p_1p_2\cdots p_m}=x^{p_1p_2\cdots p_m}+(p_1p_2\cdots p_m)n(r_1r_2\cdots r_m),$$ 
with $r_1, r_2,\dots, r_m\in R$. In other words, 
$$y^{s}=x^s+sh,$$ 
where $h=nr_1r_2\cdots r_m\in N$. Finally, since $h\in N$ and $sN=0$, it follows that $y^{s}=x^s.$ Hence, the function $H$ is well defined. It is easily verified that the function $H$ is multiplicative.
\item Since  $(g+N)^2=a+N$, it follows that $(H(g+N))^2=H(a+N)$, thus $\left(g^s\right)^2=a^s$, i.e., $a^s$ is a quadratic residue in the ring $R$.
\end{enumerate}
\cqd

An example illustrating the previous proposition is presented which it allows to discuss additional properties of the function $H$. Consider $R=\mathbb{Z}_{25}$ and $N=\langle 5 \rangle=\{0,5,10,15,20\}$. It is clear that $N$ has nilpotency  index $t=2$,  $sN={0}$ for $s=5$ and
$$
\frac{\mathbb{Z}_{25}}{\langle 5 \rangle}\cong \mathbb{Z}_5.
$$

\begin{itemize}
\item By setting $x+N=\overline{x}$, it is not difficult to see that $H(\overline{0})=0$, $H(\overline{1})=1$, $H(\overline{2})=7$, $H(\overline{3})=18$ and $H(\overline{4})=24$. Since, $\overline{0}, \overline{1}$ and $\overline{4}$ are the quadratic residues in  
the ring $\mathbb{Z}_{25}/\langle 5 \rangle$, then $0,1$ and $24$ are quadratic residues in the ring $\mathbb{Z}_{25}$.
\item Since $(3+N)^2=4+N$ in $\mathbb{Z}_{25}/\langle 5 \rangle$. From proposition \ref{quadN2}, it follows that the function $\eta:\{3,8,13,18,23\}\rightarrow  \{4,9,14,19,24\}$ given by $\eta(x)=x^2 \mod(25)$ is a bijective function, in particular, all the elements in the equivalence class $4+N$ are quadratic residues in the ring $\mathbb{Z}_{25}$. Thus, the only quadratic residue in the equivalence class $4+N$ that is mapping by the function $H$ is $24$.
\item Since $H(1+N+1+N)=7$ and $H(1+N)+H(1+N)=1+1=2$, then the function $H$ is not in general a ring homomorphism.
\end{itemize}

It is to be noticed that the hypothesis in claim \ref{port} of Proposition \ref{potencia}, which requires that all prime factors of $s$ be greater or equal than the nilpotency index $t$ of the ideal $N$, restricts enormously the number of applications of that proposition. For instance, if we consider $R=\mathbb{Z}_{2^t}$ with $2\leq t$ and $N=\langle 2 \rangle$, it is clear that $N^{t}=\{0\}$  and $sN=\{0\}$ for $s=2^{t-1}$. Thus, according to claim \ref{port} of Proposition \ref{potencia}, in order to get quadratic residues in the ring $R$ by computing the quadratic residues in the ring $R/N\cong \mathbb{Z}_2$, it is necessary that $t\leq 2$, therefore $t=2$. Hence, we can only obtain quadratic residues in the ring $\mathbb{Z}_4$, which is easily done by hand. In the following lines, we show how to overcome such restrictions.\\


\begin{defin}{\cite[Definition 3.2]{liftidemp}}\label{PosLiftIde}
We say that a collection $\{N_1,..., N_k\}$ of ideals of a ring $R$ satisfies the {\it CNC-condition} if the following properties hold:
\begin{enumerate}
\item \label{chc} {\bf Chain condition:} $\{0\}=N_k\subset N_{k-1}\subset\cdots \subset N_{2}\subset N_{1}\subset R$.
\item \label{nic} {\bf Nilpotency condition:} for $i=1,2,3,\ldots,k-1$, there exists $t_i \geq 2$ such that $N_i^{t_i}\subset N_{i+1}$.
\item \label{cac} {\bf Characteristic condition:} for $i=1,2,3,\ldots,k-1$, there exists $s_i\geq 1$ such that  $s_iN_i\subset N_{i+1}$. In addition, the prime factors of  $s_i$ are greater than or equal to $t_i$.   
\end{enumerate}
The minimum number $t_i$ satisfying the nilpotency condition will be called the nilpotency index of the ideal $N_i$ in the ideal $N_{i+1}$. Similarly,  
the minimum number $s_i$ satisfying the characteristic condition will be called the characteristic of the ideal $N_i$ in the ideal $N_{i+1}$. 
\end{defin}

\noindent
The nilpotency condition and the characteristic condition of the previous definition can be stated as follows:

\begin{itemize}
\item[a.]\label{pronic} The nilpotency condition is equivalent to the following condition: for $i=1,2,\ldots,k-1$, $N_{i}/N_{i+1}$ is a nilpotent ideal of index $t_i$ in the ring $R/N_{i+1}$, (for details see {\cite[Definition 3.2]{liftidemp}}).  
\item[b.]\label{procac} The characteristic condition is equivalent to the following condition: for $i=1,2,\ldots,k-1$, there exists a natural number $s_i\geq1$ such that $s_i(N_{i}/N_{i+1})=0$ in the ring $R/N_{i+1}$, (for details see {\cite[Definition 3.2]{liftidemp}}).
\end{itemize}

\begin{theo}\label{IdemGeral} 
Let $	R$ be a commutative ring, $\{N_1, N_2, \ldots, N_k\}$ a collection of ideals of $R$ satisfying the CNC-condition and let $s_i$ be the characteristic of the ideal $N_i$ in the ideal $N_{i+1}$. If $a+N_1$ is a quadratic residue in $R/N_1$, then $a^{s_1s_2\cdots s_{k-1}}$ is a quadratic residue in $R$. More precisely, if $g\in R$ is such that $(g+N_1)^2=a+N_1$, then
\begin{equation}\label{quadraticgeral}
\left(g^{s_1s_2\cdots s_{k-1}}\right)^2=a^{s_1s_2\cdots s_{k-1}}.
\end{equation}
\end{theo} 

\proof Note first that since the ideals $N_i$ satisfy the chain condition given in definition \ref{PosLiftIde}, for all $i=1,2,\dots,k-1$, the following isomorphism holds
\begin{equation}\label{voltare}
R/N_{i}\cong \frac{(R/N_{i+1})}{(N_i/N_{i+1})}.
\end{equation}
In addition, since the ideals $N_i$ satisfy the nilpotency condition and characteristic condition with characteristics $s_i$ respectivelly, from claim \ref{port} of proposition \ref{potencia}, it follows that for $i=1,2,\dots,k-1$, the functions
$$
H_i:R/N_i\rightarrow R/N_{i+1}, \hspace{0.5cm} H_i(x+N_i)=x^{s_i}+N_{i+1}
$$
are well defined and multiplicative. Hence, if $(g+N_1)^2=a+N_1$,
$$
H_{k-1}\circ\cdots\circ H_1((g+N_1)^2)=H_{k-1}\circ\cdots\circ H_1(a+N_1),
$$ 
whence the identity in (\ref{quadraticgeral}) is obtained.
\cqd

\begin{obs}
It follows from the proof of the Theorem (\ref{IdemGeral}) that, if $a\in R$ is such that $a+N_1$ is a quadratic residue in $R/N_1$, then $H_1(a+N_1)=a^{s_1}+N_{2}$ is a quadratic residue in $R/N_{2}$. In the same way, $H_2(a^{s_1}+N_2)=a^{s_1s_2}+N_{3}$ is a quadratic residue in $R/N_{3}$,  and so on. At the end of this process, it is obtained that $a^{s_1s_2\cdots s_{k-1}}$ is a quadratic residue in $R$. The following chain of multiplicative functions, 
$$
\frac{R}{N_1}\xrightarrow{H_1}\frac{R}{N_{2}}\xrightarrow{H_{2}} \cdots \xrightarrow{H_{k-2}}\frac{R}{N_{k-1}}\xrightarrow{H_{k-1}}\frac{R}{N_{k}}=R,\hspace{0.3cm}\text{with}\hspace{0.3cm}H_i(x+N_i)=x^{s_i}+N_{i+1}
$$
appears naturally in that process. 
\end{obs}

\begin{theo}\label{IdemGeralSeg} 
Let $	R$ be a commutative ring with identity, $\{N_1, N_2, \ldots, N_k\}$ a collection of ideals of $R$ satisfying both the Chain condition and the Nilpotency condition. Assuming that $2+N_1\in (R/N_1)^*$, the following claims hold
\begin{enumerate}
\item[1.]  $a+N_1\in q((R/N_1)^*)$ if and only if $a+N_1\subset q(R^*)$.
\item[2.] The cardinality of the set $q(R^*)$ is given by
\begin{equation}\label{card4}
|q(R^*)|=|N_1||q((R/N_1)^*)|
\end{equation}

\item[3.] If $a+N_1\in q((R/N_1)^*)$, then 
\begin{equation}\label{numsol}
s(a)=s(a+N_{k-1})=\cdots=s(a+N_1).
\end{equation}	
\item[4.] If for each $i=1,2,3,\dots,k-1$, there exists $\alpha_i$ such that, $|s(a+N_i)|=\alpha_i$ for all $a+N_i\in q((R/N_i)^*)$, then
\begin{equation}\label{nubsol}
|(R/N_{i+1})^*|=\alpha_i|q((R/N_{i+1})^*)|.
\end{equation}	
In particular,
\begin{equation}\label{nubsolta}
|R^*|=\alpha_{k-1}|q(R^*)|.
\end{equation}	

\end{enumerate}
\end{theo} 
\proof 1. It is easy to see that if $a+N_1\subset q(R^*)$, then $a+N_1\in q((R/N_1)^*)$. Now, we proceed to prove the other implication of the statement. From the isomorphism given in (\ref{voltare}), the fact that $N_i/N_{i+1}$ is a nilpotent ideal of index $t_i$ in the ring $R/N_{i+1}$ and the fact that $2+N_1\in (R/N_1)^*$, from proposition (\ref{jacobson}), it follows that 
$$
(2+N_{i+1})+N_i/N_{i+1}\in ((R/N_{i+1})/(N_i/N_{i+1}))^*
$$ for all $i\in{1,2,3,\dots,k}$. Now, let $b\in a+N_1$, since $b+N_1=a+N_1\in q((R/N_1)^*)$, it follows from the isomorphism 
\begin{equation*}\label{iso1}
R/N_1\cong \frac{(R/N_2)}{(N_1/N_2)},
\end{equation*}
that $(b+N_2)+N_1/N_2\in q(((R/N_2)/(N_1/N_2))^*)$, thus, from claim \ref{ide1} of proposition \ref{coroe}, it follows that
$$
(b+N_{2})+N_1/N_{2}=\{b+n+N_{2}| n\in N_1\} \subset q((R/N_{2})^*),
$$
in particular, it is concluded that $b+N_2\in q((R/N_{2})^*)$. Similarly, from the isomorphism
\begin{equation*}\label{iso2}
R/N_2\cong \frac{(R/N_3)}{(N_2/N_3)},
\end{equation*}
it follows that $(b+N_3)+N_2/N_3\in q(((R/N_3)/(N_2/N_3))^*)$, thus, from item \ref{ide1} of proposition \ref{coroe}, it follows that
$$
(b+N_{3})+N_2/N_{3}=\{b+n+N_{3}| n\in N_2\} \subset q((R/N_{3})^*),
$$
in particular, it is concluded that $b+N_3\in q((R/N_{3})^*)$. Continuing this process, it is finally shown that $b+N_k=\{b\}\in q((R/N_{k})^*)$, which immediately implies that $b\in q(R^*)$. This shows that $a+N_1\subset q(R^*)$, as we wanted to prove.\\
2. From the isomorphism given in (\ref{voltare}) and item \ref{ide3} of proposition \ref{coroe}, it follows that
\begin{equation}\label{estr}
|N_i/N_{i+1}||q((R/N_i)^*)|=|q(((R/N_{i+1})/(N_i/N_{i+1}))^*)|=|q((R/N_{i+1})^*)|,
\end{equation}
thus from Lagrange's theorem,
$$
|q(R^*)|=|q((R/N_{k})^*)|=\frac{|N_{k-1}|}{|N_{k}|}|q((R/N_{k-1})^*)|= \frac{|N_{k-1}|}{|N_{k}|}\frac{|N_{k-2}|}{|N_{k-1}|} \cdots \frac{|N_{1}|}{|N_{2}|}q((R/N_{1})^*), 
$$
whence the identity in (\ref{card4}) is obtained.\\
3. Again, from the isomorphism given in (\ref{voltare}), it follows that
$$
s(a+N_i)=s(a+N_{i+1}+N_i/N_{i+1})
$$
for all $i=1,2,3,\dots,k-1$. On the other hand, since $N_i/N_{i+1}$ is a nilpotent ideal of index $t_i$ in the ring $R/N_{i+1}$ and the fact that $(2+N_{i+1})+N_i/N_{i+1}\in ((R/N_{i+1})/(N_i/N_{i+1}))^*$, from claim \ref{ide2} of proposition \ref{coroe}, it follows that
$$
s(a+N_{i+1}+N_i/N_{i+1})=s(a+N_{i+1}).
$$
for $i=1,2,3,\dots,k-1$. From the previous identities, it follows that $s(a+N_i)=s(a+N_{i+1})$ for $i=1,2,3,\dots,k-1$, this of course implies the equalities appearing in (\ref{numsol}).\\
4. It follows from the isomorphism given in (\ref{voltare}) and claim \ref{ide3.1} of proposition \ref{coroe} that
$$
|(R/N_{i+1})^*|=|N_i/N_{i+1}|\sum_{a+N_i\in q((R/N_i)^*)} |s(a+N_i)|.
$$
Since, $|s(a+N_i)|=\alpha_i$, it is deduced from the former identity that
$$
|(R/N_{i+1})^*|=\alpha_i |N_i/N_{i+1}||q((R/N_i)^*)|.
$$
Finally, identity in (\ref{nubsol}) follows from (\ref{estr}).
\cqd


\section{Applications of the main results}

\quad In this section Theorems \ref{IdemGeral} and \ref{IdemGeralSeg} will be used in order to describe properties of the set of invertible quadratic residues for several classes of rings which include: rings containing a nilpotent ideal; group rings $RG$ where $R$ is a commutative ring containing a collection of ideals satisfying the CNC-condition and $G$ is a commutative group; polynomial ring $R[x]$ where $R$ is a commutative ring containing a collection of ideals satisfying the CNC-condition. Examples are given illustrating the results.

\subsection{Rings containing a nilpotent ideal}

If $R$ is a commutative ring containing a nilpotent ideal $N$, by invoking Theorems \ref{IdemGeral} and \ref{IdemGeralSeg}, properties of the set of invertible quadratic residues of the ring $R$ are described.
\begin{propo}\label{GeralNil} 
Let $R$ be a commutative ring and $N$ a nilpotent ideal of nilpotency index $k\geq 2$ in $R$. Then, the following statements hold,
\begin{enumerate}
\item\label{E1GN} Let $s>1$ be the characteristic of the quotient ring $R/N$. If $a+N$ is a quadratic residue in $R/N$, then $a^{s^{k-1}}$ is a quadratic residue in $R$. More precisely, if $g\in R$ is such that $(g+N)^2=a+N$, then
\begin{equation}\label{ticgeralnil}
\left(g^{s^{k-1}}\right)^2=a^{s^{k-1}}.
\end{equation}
\item\label{E2GN} Assuming that $2+N\in (R/N)^*$, the following claims hold,
\begin{enumerate}
\item[a).]  $a+N\in q((R/N)^*)$ if and only if $a+N\subset q(R^*)$.
\item[b).] The cardinality of the set $q(R^*)$ is given by
\begin{equation}\label{ard4nil}
|q(R^*)|=|N||q((R/N)^*)|
\end{equation}
\item[c).] If $a+N\in q((R/N)^*)$, then 
\begin{equation}\label{umsolnil}
s(a)=s(a+N^{k-1})=\cdots=s(a+N).
\end{equation}	
\item[d).] If there exists $\beta$ such that, $|s(a+N)|=\beta$ for all $a+N\in q((R/N)^*)$, then
\begin{equation}\label{ubsolnil}
|(R/N^{i+1})^*|=\beta|q((R/N^{i+1})^*)|.
\end{equation}	
In particular,
\begin{equation}\label{ubsoltanil}
|R^*|=\beta |q(R^*)|.
\end{equation}	
\end{enumerate}
\end{enumerate}
\end{propo} 
\proof In  \cite{liftunits} (Proposition 4.1), it is proven that the collection  
$$
B=\{N, N^2,...,N^k\}
$$ 
of ideals of the ring $R$ satisfies the CNC-condition with nilpotency index and characteristic of the ideal $N^{i}$ in the ideal $N^{i+1}$ being $t_i=2$ and $s_i=s$ for all $i=1,2,3,\dots,k-1$. Therefore, the proof of this proposition is now a clear consequence of Theorems \ref{IdemGeral} and \ref{IdemGeralSeg}. 
\cqd

\begin{example}\label{Horacio}\rm{
Let $p$ be an odd prime number, $i\in\mathbb{N}$ and let $R=\{a+bu : a,b\in \mathbb{Z}_{p^i}, u^2=0 \}$. It is readily seen that $R$ with the (obvious) addition and multiplication operations is a commutative ring with cardinality $|R|=p^{2i}$. It is also easily seen that $R$ is isomorphic to the ring of polynomials with coefficients in  $\mathbb{Z}_{p^i}$ modulo the ideal generated by $x^2$, that is $\mathbb{Z}_{p^i}[x]/\langle x^2\rangle$. It is readely seen that
$$
R^*=\{a+bu : a\in(\mathbb{Z}_{p^i})^*, b\in \mathbb{Z}_{p^i}\},
$$
so the cardinality of $R^*$ is $|R^*|=\varphi(p^i)p^i=(p-1)p^{2i-1}$, where $\varphi$ denotes the Euler totient function. On the other hand, it is verified that the ideal $N=\langle p,u\rangle$ has nilpotency index $k=i+1$ and that $|N|=p^{2i-1}$, then it follows that $N$ is a maximal ideal of $R$ with
$$
\frac{R}{N}\cong \mathbb{Z}_p,
$$  
whence $|(R/N)^*|=p-1$ and the characteristic of the quotient ring $R/N$ is $s=p$. From the latter isomorphism and proposition \ref{GeralNil}, it is concluded that 
\begin{itemize}
\item $a+bu\in q(R^*)$ if and only if $a\mod(p)\in q(\mathbb{Z}_p^*)$. 
\item Let $a+bu\in R$ if $a\mod(p)\in q(\mathbb{Z}_p^*)$ then for all $b\in \mathbb{Z}_{p^i}$
$$
(a+bu)^{p^{i}}=(a\mod (p))^{p^{i}}
$$
is an invertible quadratic residue in $R$.
\item The number of invertible quadratic residues of the ring $R$ is given by
$$
|q(R^*)|=|N||q((R/N)^*)|=\frac{p^{2i-1}(p-1)}{2}
$$
\item Let $a+bu\in R$, if $a\mod(p)\in q(\mathbb{Z}_p^*)$ then for all $b\in \mathbb{Z}_{p^i}$ the number of solutions in $R$ of the equation $x^2=a+bu$ is equal to 2, in other words
$$
s(a+bu)=2.
$$
\end{itemize}
 }
\end{example}

An easy application of the previous result is the following:

\begin{coro}\label{coroeta}
Let $R$ be a commutative ring and $c$ a nilpotent element of index $k\geq 2$ in $R$. Then, the following statements hold:
\begin{enumerate}
\item\label{E1GNcor} Let $s>1$ be the characteristic of the quotient ring $R/\langle c\rangle$. If $a+\langle c\rangle$ is a quadratic residue in $R/\langle c\rangle$, then $a^{s^{k-1}}$ is a quadratic residue in $R$. More precisely, if $g\in R$ is such that $(g+\langle c\rangle)^2=a+\langle c\rangle$, then
\begin{equation}\label{icgeralnic}
\left(g^{s^{k-1}}\right)^2=a^{s^{k-1}}.
\end{equation}
\item\label{E2GNcor} Assuming that $2+\langle c\rangle\in (R/\langle c\rangle)^*$, the following claims hold:
\begin{enumerate}
\item[a.]  $a+\langle c\rangle\in q((R/\langle c\rangle)^*)$ if and only if $a+\langle c\rangle\subset q(R^*)$.
\item[b.] The cardinality of the set $q(R^*)$ is given by
\begin{equation}\label{rd4c}
|q(R^*)|=|\langle c\rangle||q((R/\langle c\rangle)^*)|.
\end{equation}
\item[c.] If $a+\langle c\rangle\in q((R/\langle c\rangle)^*)$, then 
\begin{equation}\label{msolc}
s(a)=s(a+\langle c^{k-1}\rangle)=\cdots=s(a+\langle c\rangle).
\end{equation}	
\item[d.] If there exists $\beta$ such that $|s(a+\langle c\rangle)|=\beta$ for all $a+\langle c\rangle\in q((R/\langle c\rangle)^*)$, then
\begin{equation}\label{bsolc}
|(R/\langle c^{i+1}\rangle)^*|=\beta|q((R/\langle c ^{i+1}\rangle)^*)|.
\end{equation}	
In particular,
\begin{equation}\label{bsoltac}
|R^*|=\beta |q(R^*)|.
\end{equation}	
\end{enumerate}
\end{enumerate}
\end{coro} 
\proof Since $R$ is a commutative ring, $\langle c \rangle$ is a nilpotent ideal of nilpotency index $k$ in $R$, and the result follows immediately from Proposition \ref{GeralNil}
\medskip
\cqd


\subsection{Group rings}

If $R$ is a commutative ring containing a collection of ideals satisfying the CNC-condition and $G$ is a commutative group, by invoking Theorems \ref{IdemGeral} and \ref{IdemGeralSeg}, properties of the set of invertible quadratic residues of the group ring $RG$ are described. 

\begin{propo}\label{GeralGR} 
Let $	R$ be a commutative ring and $G$ a commutative group. Let $\{N_1, N_2, \ldots, N_k\}$ be a collection of ideals of $R$ satisfying the CNC-condition. Then, the following statements hold:
\begin{enumerate}
\item\label{E1GR} Let $s_i$ be the characteristic of the ideal $N_i$ in the ideal $N_{i+1}$. If $a+N_1G$ is a quadratic residue in $(R/N_1)G$, then $a^{s_1s_2\cdots s_{k-1}}$ is a quadratic residue in $RG$. More precisely, if $g\in RG$ is such that $(g+N_1G)^2=a+N_1G$, then
\begin{equation}\label{eralnigr}
\left(g^{s_1s_2\cdots s_{k-1}}\right)^2=a^{s_1s_2\cdots s_{k-1}}.
\end{equation}
\item\label{E2GR} Assuming that $2+N_1G\in ((R/N_1)G)^*$, then the following claims hold:
\begin{enumerate}
\item[a.]  $a+N_1G\in q(((R/N_1)G)^*)$ if and only if $a+N_1G\subset q((RG)^*)$.
\item[b.] The cardinality of the set $q((RG)^*)$ is given by
\begin{equation}\label{rd4gr}
|q((RG)^*)|=|N_1|^{|G|}|q(((R/N_1)G)^*)|.
\end{equation}
\item[c.] If $a+N_1G\in q(((R/N_1)G)^*)$ then, 
\begin{equation}\label{msolgr}
s(a)=s(a+N_{k-1}G)=\cdots=s(a+N_1G).
\end{equation}	
\item[d.] If there exists $\beta$ such that $|s(a+N_1G)|=\beta$ for all $a+N_1G\in q(((R/N_1)G)^*)$ then,
\begin{equation}\label{bsolgr}
|((R/N_{i+1})G)^*|=\beta|q(((R/N_{i+1})G)^*)|,
\end{equation}	
for $i=1,2,\cdots,k-1$. In particular,
\begin{equation}\label{bsoltagr}
|(RG)^*|=\beta |q((RG)^*)|.
\end{equation}	
\end{enumerate}
\end{enumerate}
\end{propo} 

\proof 
In  \cite{liftunits} (Proposition 4.9), it is shown that the collection  
$$
B=\{N_1G, N_2G, \ldots, N_kG\}
$$ 
of ideals of the ring $RG$ satisfies the CNC-condition with nilpotency index and characteristic of the ideal $N_iG$ in the ideal $N_{i+1}G$ being exactly the same nilpotency index and characteristic of the ideal $N_i$ in the ideal $N_{i+1}$. 
Therefore, the proof of this proposition is a direct conequence of Theorems \ref{IdemGeral} and \ref{IdemGeralSeg}.
\cqd

\begin{coro}\label{GeralNilGR} 
Let $G$ be a commutative group, $R$ be a commutative ring and $N$ a nilpotent ideal of index $k$ in $R$. Then, the following statements hold:
\begin{enumerate}
\item\label{E1GRNC} Let $s>1$ be the characteristic of the quotient ring $R/N$. If $a+NG$ is a quadratic residue in $(R/N)G$, then $a^{s^{k-1}}$ is a quadratic residue in $RG$. More precisely, if $g\in RG$ is such that $(g+NG)^2=a+NG$, then
\begin{equation}\label{ralnigrc}
\left(g^{s^{k-1}}\right)^2=a^{s^{k-1}}.
\end{equation}
\item\label{E2GRNC} Assuming that $2+NG\in ((R/N)G)^*$, the following claims hold
\begin{enumerate}
\item[a.]  $a+NG\in q(((R/N)G)^*)$ if and only if $a+NG\subset q((RG)^*)$.
\item[b.] The cardinality of the set $q((RG)^*)$ is given by
\begin{equation}\label{d4grc}
|q((RG)^*)|=|N|^{|G|}|q(((R/N)G)^*)|.
\end{equation}
\item[c.] If $a+NG\in q(((R/N)G)^*)$, then 
\begin{equation}\label{solgrc}
s(a)=s(a+N^{k-1}G)=\cdots=s(a+NG).
\end{equation}	
\item[d.] If there exists $\beta$ such that $|s(a+NG)|=\beta$ for all $a+NG\in q(((R/N)G)^*)$ then,
\begin{equation}\label{solgrc}
|((R/N^{i+1})G)^*|=\beta|q(((R/N^{i+1})G)^*)|,
\end{equation}	
for $i=1,2,\cdots,k-1$. In particular,
\begin{equation}\label{soltagrc}
|(RG)^*|=\beta |q((RG)^*)|.
\end{equation}	
\end{enumerate}
\end{enumerate}

\end{coro}
\proof 
The proof of this corollary is a direct consequence of Proposition \ref{GeralGR} and the fact that the collection $\{N, N^2,...,N^{k}\}$ of ideals of the ring $R$ satisfies the CNC-condition with constant characteristic $s_i=s$ for all $i=1,2,3,\cdots, k-1$.
\cqd 

\begin{example}\label{destricta}\rm{
Let $p$ be an odd prime number, $i\in\mathbb{N}$ and let $R=\{a+bu : a,b\in \mathbb{Z}_{p^i}, u^2=1 \}$ be the group ring $\mathbb{Z}_{p^i}G$ where $G=\{1,u\}$ is the cyclic group of order $n=2$. It is readily seen that $R$ with the (obvious) addition and multiplication operations is a commutative ring with cardinality $|R|=p^{2i}$. It is also easily seen that $R$ is isomorphic to the ring of polynomials with coefficients in  $\mathbb{Z}_{p^i}$ modulo the ideal generated by $x^2-1$ in $R$, that is $\mathbb{Z}_{p^i}G\cong \mathbb{Z}_{p^i}[x]/\langle x^2-1\rangle$. It is readily seen that
$$
(\mathbb{Z}_pG)^*=\{a+bu : a\neq b, a\neq -b\},
$$
so the cardinality of $(\mathbb{Z}_pG)^*$ is $|(\mathbb{Z}_pG)^*|=(p-1)^2$. In addition, since $N=\langle p \rangle$ has nilpotency index $k=i$ in $\mathbb{Z}_{p^i}$, and 
$$
\frac{\mathbb{Z}_{p^i}G}{\langle p\rangle G}\cong \mathbb{Z}_pG,
$$
then it is deduced that $|(\mathbb{Z}_{p^i}G)^*|=|(\mathbb{Z}_pG)^*||\langle p \rangle G|=(p-1)^2p^{2(i-1)}$. From the latter isomorphism and the proposition \ref{GeralNilGR}, it is concluded that: 
\begin{itemize}
\item $a+bu\in q((\mathbb{Z}_{p^i}G)^*)$ if and only if $(a\mod(p))+(b\mod(p))u\in q((\mathbb{Z}_{p}G)^*)$. 
\item Let $a+bu\in \mathbb{Z}_{p^i}G$, if  $(a\mod(p))+(b\mod(p))u\in q((\mathbb{Z}_{p}G)^*)$, then 
$$
(a+bu)^{p^{i-1}}=((a\mod(p))+(b\mod(p))u)^{p^{i-1}}
$$
is an invertible quadratic residue in $R=\mathbb{Z}_{p^i}G$.
\item The number of invertible quadratic residues of the ring $R$ is given by
$$
|q((\mathbb{Z}_{p^i}G)^*)|=|N|^{|G|}|q((\mathbb{Z}_pG)^*)|=p^{2(i-1)}|q((\mathbb{Z}_pG)^*)|
$$
\item Let $a+bu\in \mathbb{Z}_{p^i}G$, if $(a\mod(p))+(b\mod(p))u\in q((\mathbb{Z}_{p}G)^*)$ and $|s((a\mod(p))+(b\mod(p))u)|=\beta$, then 
$$
|s(a+bu)|=\beta.
$$
If additionally, $|s((a\mod(p))+(b\mod(p))u)|=\beta$ for all $(a\mod(p))+(b\mod(p))u\in q((\mathbb{Z}_{p}G)^*)$ then,
$$
|(\mathbb{Z}_{p^i}G)^*|=\beta |q((\mathbb{Z}_{p^i}G)^*)|.
$$
\end{itemize}
For instance, if $p=3$, it is easy to see that $q((\mathbb{Z}_{3}G)^*)=\{1\}$ and the number of solutions in $\mathbb{Z}_{3}G$ of the equation $x^2=1$ is equal to $4$, in other words $|s(1)|=4$. Thus, if $a+bu\in \mathbb{Z}_{3^i}G$ is such that $a\equiv 1 \mod(3)$ and $b\equiv 0 \mod(3)$, then 
$$
(a+bu)^{3^{i-1}}=1,
$$
$ |s(a+bu)|=4$, $|q((\mathbb{Z}_{3^i}G)^*)|=3^{2(i-1)}$ and $|(\mathbb{Z}_{3^i}G)^*|=(4) 3^{2(i-1)}$.
}
\end{example}

\end{document}